\documentclass[12pt]{amsart}
\usepackage{latexsym, amsmath, amscd, amssymb, amsthm,longtable} 
\usepackage[T2A]{fontenc}
 \usepackage[cp1251]{inputenc}
\usepackage[english]{babel}
\usepackage{color}
\usepackage{longtable}
 \def\mathbi#1{\textbf{\em #1}}
 
\begin{document}

\noindent

 \title[  Hilbert polynomials] { Hilbert polynomials of the algebras of  $SL_2$-invariants }

\author{Leonid Bedratyuk}

\begin{abstract}
We offer a Maple-procedure for computing of the Hilbert polynomials of the algebras of  $SL_2$-invariants
\end{abstract}
\maketitle

{\bf 1.} Let $K$ be a field, ${\rm char}K=0.$ Let $V_d $ be  $d+1$-dimensional  $SL_2$-module of binary forms of degree $d$ and let $V_{\mathbi{d}} =
V_{d_1} \oplus V_{d_2} \oplus \cdots \oplus V_{d_n},$ ${\mathbi{d}=(d_1,d_2,\ldots, d_s).}$ Denote by  $K[V_{\mathbi{d}}]^{SL_2}$   the algebra of polynomial $SL_2$-invariant functions on $V_{\mathbi{d}}.$  It is well-known that the algebra
 $\mathcal{I}_{\mathbi{d}}:=K[V_{\mathbi{d}}]^{SL_2}$ is finitely generated and  graded  : 
$$
\mathcal{I}_{\mathbi{d}}=(\mathcal{I}_{\mathbi{d}})_0 \oplus (\mathcal{I}_{\mathbi{d}})_1 \oplus (\mathcal{I}_{\mathbi{d}})_2 \oplus \cdots \oplus (\mathcal{I}_{\mathbi{d}})_n \oplus +\cdots,
$$
where  $(\mathcal{I}_{\mathbi{d}})_n$  is a vector $K$-space of invariants of  degree $n.$ The Hilbert function of the algebra    $\mathcal{I}_{\mathbi{d}}$ is  defined as dimension of the vector space  $(\mathcal{I}_{\mathbi{d}})_n:$
$$
\mathcal{H}(\mathcal{I}_{\mathbi{d}},n)=\dim (\mathcal{I}_{\mathbi{d}})_n.
$$

It is well-known   {\cite{{St-30}}-\cite{Eis_Syz}}  that the Hilbert function  of finitely generated graded $K$-algebra is equal (starting from some $n$) to a  polynomial of $n:$   
$$
\mathcal{H}(\mathcal{I}_{\mathbi{d}},n)=h_0(n)n^{r}+h_1(n)n^{r-1} +\cdots, 
$$
where $h_i(n)$  is some  periodic function  with values in $\mathbb{Q}.$  Then such a polynomial is  called the  Hilbert polynomial  of  graded algebra. From combinatorial point of view  the Hilbert polynomials are so-called  quasi-polynomials, see   \cite{St-En_1}, Chapter 4.

 In this short  notes we present a Maple-procedure  for  calculation  of    the Hilbert polynomials   for  the algebra $\mathcal{I}_{\mathbi{d}}.$    

{\bf 2.} For  the case    of  one binary form ($s=1$)  there exists  (\cite{Hilb}, \cite{SP}) classical Cayley-Sylvester formula for  calculation of values of Hilbert function of  $\mathcal{I}_{d}:$
$$
\mathcal{H}(\mathcal{I}_{d},n)=\omega_d \left(n,0\right)- \omega_d \left(n,2\right),
$$ 
where $\omega_d(n,k) $  is the number non-negative integer solutions of the system: $$
\left \{ \begin{array}{l} \alpha_1+2\alpha_2+\cdots + d\, \alpha_d=\displaystyle \frac{d\,n-k}{2},\\ \alpha_0+\alpha_1+\cdots +  \alpha_d=n.
\end{array}
\right.
$$  
For generalisation of the Cayley-Sylvester formula to $\mathcal{I}_{\mathbi{d}}$  see  \cite{B-SC}.
Also, see \cite{Hilb}, \cite{SP}, we  have:
$$
\mathcal{H}(\mathcal{I}_{d},n)=\left[q^{\frac{n\,d}{2}}\right] \left((1-q)\left[ \begin{array}{c} d \\ n  \end{array} \right]_q\right),
$$
here  $\left[ \begin{array}{c} d \\ n  \end{array} \right]_q$ is $q$-binomial coefficient:
$$
\left[ \begin{array}{c} d \\ n  \end{array} \right]_q:=\frac{(1-q^{d+1})(1-q^{d+2})\ldots(1-q^{d+n})}{(1-q)(1-q^{2})\ldots(1-q^{n})},
$$
and  $\left[q^{\frac{n\,d}{2}}\right]$ denotes the coefficient of  $q^{\frac{n\,d}{2}}.$  
But  it  seems  impossible to extract    the Hilbert polynomials from this representation of the Hilbert function.

To calculate the Hilbert polynomials we have used our knowledge of Poincare series of the algebras of invariants  $\mathcal{I}_{\mathbi{d}},$ see \cite{B-Poincare}.
The  following statement is  our  main computational tool :

\noindent
{\bf Theorem.}{\it \,  Let the Poincar\'e series of a finitely generated graded $K$-algebra  $A$   has the form
\begin{gather*}
\mathcal{P}(A,z)=\frac{R(z)}{Q(z)},
\end{gather*}
the polynomials  $R(z),Q(z),$ $\deg R(z) < \deg Q(z)$   are coprime and let  $\lambda_1, \lambda_2, \ldots \lambda_r$ are the roots of the denominator $Q(z)$ with the multiplicities  $k_1,k_2, \ldots, k_r.$ Then the Hilbert polynomials of $A$ has the form
$$
\mathcal{H}(A,n)=\tau_{k_1}(n)\lambda_1^{k_1}+\tau_{k_2}(n)\lambda_1^{k_2}+\cdots+\tau_{k_r}(n)\lambda_1^{k_r}, n=0,1,2,\ldots, 
$$
here  $\tau_{k_i}(n)$ is a polynomial  of  $n$ with rational coefficients of degree  $ \leq k_i-1.$
}

 The proof follows from  \cite{St-En_1} Theorem 4.1.1  and  Proposition  4.1.1

\noindent
{\bf Example.} Let us calculate the Hilbert  polynomial for the algebra of  joint  invariants of  two binary forms of degrees 2 and 3. The  Poincare series  equals:
\begin{gather*}
\mathcal{P}(\mathcal{I}_{2,3},z)={\frac {{z}^{6}-{z}^{5}+{z}^{4}-{z}^{3}+{z}^{2}-z+1}{\left( 1-{z}^{3} \right)  \left( 1-{z}^{5} \right)  \left(1-z+
{z}^{2}-z^3 \right)  \left( 1-{z}^{2} \right)}}=\\=1+{z}^{2}+{z}^{3}+2\,{z}^{4}+2\,{z}^{5}+3\,{z}^{6}+4\,{z}^{7}+5\,{z}^{
8}+6\,{z}^{9}+8\,{z}^{10}+9\,{z}^{11}+12\,{z}^{12}+\cdots
\end{gather*}
The  denominator 
$$
Q(z)= \left( 1-{z}^{3} \right)  \left( 1-{z}^{5} \right)  \left(1-z+
{z}^{2}-z^3 \right)  \left( 1-{z}^{2} \right) ={z}^{13}-{z}^{12}-{z}^{10}+{z}^{7}+{z}^{6}-{z}^{3}-z+1,
$$
has the root  $1$  with  the multiplicity $4$ (the transendence degree of $\mathcal{I}_{2,3}$) and the 9 roots 
\begin{gather*}
-1,i,-i,-\frac 1 4+\frac 1 4\,\sqrt {5}+\frac 1 4\,i\sqrt {10+2\,\sqrt {5}},-\frac 1 4-\frac 1 4\,\sqrt {5}+\frac 1 4\,i\sqrt {10-2\,\sqrt {5}},\\ -\frac 1 4-\frac 1 4\,\sqrt {5}-\frac 1 4\,i
\sqrt {10-2\,\sqrt {5}},-\frac 1 4+\frac 1 4\,\sqrt {5}-\frac 1 4\,i\sqrt {10+2\,\sqrt {
5}},-\frac 1 2+\frac 1 2\,i\sqrt {3},-\frac 1 2-\frac 1 2\,i\sqrt {3},
\end{gather*}
each of  multiplicity $1.$  Thus, we  are looking the Hilbert polynomial in the following form
\begin{gather*}
\mathcal{H}(\mathcal{I}_{2,3},n)=C_0+C_1 n+C_2 n^2+C_3 n^3+C_4 +C_{{5}} \left( -\frac 1 4-\frac{\sqrt {5}}{ 4}-\frac 1 4\,i\sqrt {10-2\,\sqrt {5}}
 \right) +\\+C_{{6}} \left( -\frac 1 4-\frac 1 4\,\sqrt {5}+\frac 1 4\,i\sqrt {10-2\,\sqrt 
{5}} \right) +C_{{7}} \left( -\frac 1 4+\frac 1 4\,\sqrt {5}-\frac 1 4\,i\sqrt {10+2\,
\sqrt {5}} \right) -iC_{{8}}+\\+C_{{9}} \left( -\frac 1 2+\frac 1 2\,i\sqrt {3}
 \right) +C_{{10}} \left( -\frac 1 2-\frac 1 2\,i\sqrt {3} \right) +iC_{{11}}+C_{{
12}} \left( -\frac 1 4+\frac{\sqrt {5}}{4}\,+\frac 1 4\,i\sqrt {10+2\,\sqrt {5}} \right).
\end{gather*} 
Taking into account the initial conditions $\mathcal{H}(\mathcal{I}_{2,3},0)=1,$ $\mathcal{H}(\mathcal{I}_{2,3},1)=0,$ $\mathcal{H}(\mathcal{I}_{2,3},2)=1,  $ $\ldots$ $\mathcal{H}(\mathcal{I}_{2,3},12)=12$  we  solve the linear system of equations  for $C_i$ and after simplification get the Hilbert polynomial:

\begin{longtable}{l}
$\displaystyle \mathcal{H}(\mathcal{I}_{2,3},n)={\frac {1}{360}}\,{n}^{3}+{\frac {7}{240}}\,{n}^{2}+\frac 1 6\,n+\frac{1}{10}\,\cos \left( 
\frac 2 5\,\pi \,n \right) -{\frac {1}{50}}\,\sqrt {5}\cos \left( \frac 2 5\,\pi 
\,n \right) $ +
 \\ +$\displaystyle {\frac {1}{50}}\,\sqrt {50+10\,\sqrt {5}}\sin
 \left( \frac{2\,\pi \,n}{5} \right) +\frac{1}{10}\,\cos \left( \frac{4\,\pi \,n}{5} \right) +
{\frac {\sqrt {5}}{50}}\cos \left( \frac{ 4\,\pi \,n}{ 5} \right)+\frac 1
 8\,\cos \left( \frac{ 1 \,\pi \,n}{ 2}\right){-} $ 
 \\   $\displaystyle -\frac{1}{25}\,\sqrt {10-2\,\sqrt {5}}\sin \left( \frac 4 5\,\pi \,n
 \right) +\frac 1 9\,\cos \left( \frac 2 3\,\pi \,n \right) +\frac{1}{27}\,\sqrt {3}\sin
 \left( \frac 2 3\,\pi \,n \right) -\frac 1 8\,\sin \left( \frac 1 2\,\pi \,n \right)+ $
 \\ \\  $\displaystyle  +{\frac {7}{32}}\,\cos \left( 
\pi \,n \right) -\frac{1}{25}\,\sqrt {10+2\,\sqrt {5}}\sin \left( \frac{2}{5}\,\pi \,n
 \right)-{\frac {1}
{50}}\,\sqrt {5}\sqrt {10-2\,\sqrt {5}}\sin \left( \frac 4 5\,\pi \,n
 \right)+{\frac {497}{1440}}.$
\end{longtable}

{\bf 3.} To calculate  the Hilbert  polynomial of the $SL_2$-module $V_{\mathbi{d}},$ ${\mathbi{d}=(d_1,d_2,\ldots, d_n)}$ we offer  the procedure {\tt HilbertPol} that included in our  updated Maple-package   {\tt Poincare\_Series}, \cite{B-Poincare}: 

{\noindent
{\em Command name}: {\tt HilbertPol}\\
\noindent{\em Feature}: Computes the Hilbert polynomial for the algebras of   joint invariants for the binary forms of degrees $d_1, d_2, \ldots,d_n.$\\
\noindent{\em Calling sequence}:
{\tt
HilbertPol$\left([d_1, d_2, \ldots,d_n] \right);$
}\\
\noindent{\em Parameters}:\\
{\raggedright
\begin{tabular}{lcl}
$ [d_1, d_2, \ldots,d_n]$ & - & a list of degrees of $n$  binary forms.\\
$ n $ & - & an integer, $n\geq 1.$\\

\end{tabular}\\
Below are listed of results of some calculations.

\vspace{5mm}
3.1 $\mathcal{H}(\mathcal{I}_{2},n)$

{\tt > \textcolor{red}{dd:=[2]:HilbertPol(dd);}}\\
$$\textcolor{blue}{
\frac 1 2\,\cos \left( \pi \,n \right) +\frac 12} $$

\vspace{5mm}
3.2  $\mathcal{H}(\mathcal{I}_{3},n)$

{\tt > \textcolor{red}{dd:=[3]:HilbertPol(dd);}}\\
$$\textcolor{blue}{
\frac  14\,\cos \left( \pi \,n \right) +\frac 12\,\cos \left( \frac 12\,\pi \,n
 \right) +\frac 14
} $$

\vspace{5mm}
3.3  $\mathcal{H}(\mathcal{I}_{4},n)$

{\tt > \textcolor{red}{dd:=[4]:HilbertPol(dd);}}\\
$$\textcolor{blue}{
\frac 16\,n+\frac 14\,\cos \left( \pi \,n \right) +\frac 13\,\cos \left( \frac23\,\pi \,n
 \right) -\frac 19\,\sqrt {3}\sin \left( \frac 23\,\pi \,n \right) +{
\frac {5}{12}}
} $$

\vspace{5mm}
3.4 $\mathcal{H}(\mathcal{I}_{5},n)$

{\tt > \textcolor{red}{dd:=[5]:HilbertPol(dd);}}\\

\textcolor{blue}{ \begin{gather*} 
\left( {\frac {1}{384}}\,\cos \left( \pi \,n \right) +{\frac {1}{384}} \right) {n}^{2}+\left( {\frac {1}{64}}\,\cos \left( \pi \,n \right) +{\frac {1}{64}}
 \right) n+{
\frac {3}{32}}\,\cos \left( \frac 12\,\pi \,n \right) +{\frac {107}{576}}\,\cos \left( \pi \,n \right) +\\+\frac 19
\,\cos \left( \frac 13\,\pi \,n \right) +\frac 19\,\cos \left( \frac 23\,\pi \,n
 \right) +\frac{1}{16}\,\sin \left( \frac 14\,\pi \,n \right) +\frac{1}{16}\,\cos \left( \frac 1
4\,\pi \,n \right) -\\-\frac{1}{16}\,\sin \left( \frac{3}{4}\,\pi \,n \right) +\frac{1}{16}\,\cos
 \left( \frac 34\,\pi \,n \right) +{\frac {9}{32}}\,\cos \left( \frac 12\,\pi \,
n \right) + {\frac {107}{576}}
\end{gather*}
}

\vspace{5mm}
3.5  $\mathcal{H}(\mathcal{I}_{6},n)$

{\tt > \textcolor{red}{dd:=[6]:HilbertPol(dd);}}\\

\textcolor{blue}{ \begin{gather*} 
{\frac {1}{1440}}\,{n}^{3}+ \left( {\frac {1}{64}}\,\cos \left( \pi \,
n \right) +{\frac {7}{960}} \right) {n}^{2} + \left( {\frac {7}{64}}\,\cos \left( \pi 
\,n \right) +{\frac {37}{320}} \right) n+\\+{\frac {1}{100}}\,\sqrt {5}\sqrt {10+2\,\sqrt {5}}\sin
 \left( \frac 2 5\,\pi \,n \right) -{\frac {1}{100}}\,\sqrt {10+2\,\sqrt {5}
}\sin \left( \frac 2 5\,\pi \,n \right) -{\frac {1}{50}}\,\sqrt {5}\cos
 \left( \frac 4 5\,\pi \,n \right) +\\+\frac{1}{10}\,\cos \left( \frac 45\,\pi \,n \right) -
{\frac {1}{100}}\,\sqrt {5}\sqrt {10-2\,\sqrt {5}}\sin \left( \frac 4 5\,
\pi \,n \right) -{\frac {1}{100}}\,\sqrt {10-2\,\sqrt {5}}\sin \left( 
\frac 4 5\,\pi \,n \right) +\\+\frac 1 9\,\cos \left( \frac 2 3\,\pi \,n \right) +\frac{1}{27}\,
\sqrt {3}\sin \left( \frac 2 3\,\pi \,n \right) -\frac{1}{16}\,\sin \left( \frac 1 2\,\pi 
\,n \right) +\frac{1}{16}\,\cos \left( \frac 1 2\,\pi \,n \right) +{\frac {9}{32}}\,
\cos \left( \pi \,n \right) +\\+\frac{1}{10}\,\cos \left( \frac 2 5\,
\pi \,n \right) +{\frac {1}{50}}\,\sqrt {5}\cos \left( \frac 2 5\,\pi \,n
 \right)+{\frac {497}{1440}}
\end{gather*}
}

\vspace{5mm}
3.6  $\mathcal{H}(\mathcal{I}_{(1,1)},n)$

{\tt > \textcolor{red}{dd:=[1,1]:HilbertPol(dd);}}\\
$$\textcolor{blue}{
\frac 1 2\,\cos \left( \pi \,n \right) +\frac 12} $$

\vspace{5mm}
3.7  $\mathcal{H}(\mathcal{I}_{(1,2)},n)$

{\tt > \textcolor{red}{dd:=[1,2]:HilbertPol(dd);}}\\
$$\textcolor{blue}{
\frac 16\,n+\frac 14\,\cos \left( \pi \,n \right) +\frac 13\,\cos \left( \frac23\,\pi \,n
 \right) -\frac 19\,\sqrt {3}\sin \left( \frac 23\,\pi \,n \right) +{
\frac {5}{12}}
} $$

\vspace{5mm}
3.8  $\mathcal{H}(\mathcal{I}_{(1,3)},n)$

{\tt > \textcolor{red}{dd:=[1,3]:HilbertPol(dd);}}\\
\textcolor{blue}{\begin{gather*}
\left( {\frac {1}{64}}\,\cos \left( \pi \,n \right) +{\frac {1}{64}} \right) {n}^{2}+ \left( \frac{3}{16}\,\cos
 \left( \frac 1 2\,\pi \,n \right) +{\frac {3}{32}}\,\cos \left( \pi \,n
 \right) +{\frac {3}{32}} \right) n +\\+{
\frac {7}{32}}\,\cos \left( \pi \,n \right) +{\frac {9}{16}}\,\cos \left( \frac 1 2\,\pi \,n \right)+{\frac {7}{32}}
\end{gather*}}

\vspace{5mm}
3.9  $\mathcal{H}(\mathcal{I}_{(1,4)},n)$

{\tt > \textcolor{red}{dd:=[1,4]:HilbertPol(dd);}}\\
\textcolor{blue}{\begin{gather*}
{\frac {1}{540}}\,{n}^{3}+{\frac {7}{360}}\,{n}^{2} + \left( \frac{1}{27}\,\cos \left( \frac 2 3\,
\pi \,n \right) -\frac{1}{27}\,\sqrt {3}\sin \left( \frac 2 3\,\pi \,n \right) +{
\frac {79}{540}} \right) n+\\+{\frac {1}{100}}\,\sqrt {5}\sqrt {10+2\,\sqrt {5}}
\sin \left( \frac 2 5\,\pi \,n \right) -{\frac {1}{50}}\,\sqrt {5}\cos
 \left( \frac 4 5\,\pi \,n \right) +\frac{1}{10}\,\cos \left( \frac 4 5\,\pi \,n \right) -\\-
{\frac {1}{100}}\,\sqrt {10-2\,\sqrt {5}}\sin \left( \frac 4 5\,\pi \,n
 \right) -{\frac {1}{100}}\,\sqrt {5}\sqrt {10-2\,\sqrt {5}}\sin
 \left( \frac 4 5\,\pi \,n \right) +{\frac {7}{27}}\,\cos \left( \frac 2 3\,\pi \,
n \right) -\\-{\frac {7}{81}}\,\sqrt {3}\sin \left( \frac 2 3\,\pi \,n \right) 
+\frac{3}{16}\,\cos \left( \pi \,n \right) +\frac{1}{10}\,\cos \left( 
\frac 2 5\,\pi \,n \right) +{\frac {1}{50}}\,\sqrt {5}\cos \left( \frac 2 5\,\pi 
\,n \right) -\\-{\frac {1}{100}}\,\sqrt {10+2\,\sqrt {5}}\sin \left( \frac 2 5
\,\pi \,n \right)+{\frac {763}{2160}}
\end{gather*}}

\vspace{5mm}
3.10  $\mathcal{H}(\mathcal{I}_{(1,2,3)},n)$

{\tt > \textcolor{red}{dd:=[1,2,3]:HilbertPol(dd);}}\\
\textcolor{blue}{\begin{gather*}
{\frac {17}{172800}}\,{n}^{5}+{\frac {17}{7680}}\,{n}^{4}+{\frac {1033
}{51840}}\,{n}^{3}+ \left( {\frac {1}{512}}\,\cos \left( \pi \,n
 \right) +{\frac {689}{7680}} \right) {n}^{2}+\\+ 
 \left( {\frac {4}{81}}\,\cos \left( \frac  2 3\,\pi \,n \right) -{
\frac {3}{64}}\,\sin \left( \frac  1 2\,\pi \,n \right) +{\frac {3}{64}}\,
\cos \left( \frac  1 2\,\pi \,n \right) +{\frac {9}{512}}\,\cos \left( \pi \,
n \right) +{\frac {46667}{207360}} \right) n 
 -\\-{\frac {1}{50}}\,\sqrt {
5}\cos \left( \frac  2 5\,\pi \,n \right) +{\frac {3}{50}}\,\cos \left( \frac  2 5\,
\pi \,n \right) -{\frac {1}{100}}\,\sqrt {10+2\,\sqrt {5}}\sin \left( 
\frac  2 5\,\pi \,n \right) +\\+{\frac {1}{500}}\,\sqrt {5}\sqrt {10+2\,\sqrt {5
}}\sin \left( \frac  2 5\,\pi \,n \right) +{\frac {3}{50}}\,\cos \left( \frac  4 5\,
\pi \,n \right) +{\frac {1}{50}}\,\sqrt {5}\cos \left( \frac  4 5\,\pi \,n
 \right) -\\-{\frac {1}{500}}\,\sqrt {5}\sqrt {10-2\,\sqrt {5}}\sin
 \left( \frac  4 5\,\pi \,n \right) -{\frac {1}{100}}\,\sqrt {10-2\,\sqrt {5}
}\sin \left( \frac  4 5\,\pi \,n \right) +\frac  2 9\,\cos \left( \frac  2 3\,\pi \,n
 \right) +\\+{\frac {4}{243}}\,\sqrt {3}\sin \left( 2/3\,\pi \,n \right) 
-\frac{3}{16}\,\sin \left( \frac 1 2\,\pi \,n \right) +{\frac {15}{64}}\,\cos
 \left( \frac  1 2\,\pi \,n \right) +{\frac {141}{1024}}\,\cos \left( \pi \,n
 \right) + \\+{\frac {65827}{230400}}
\end{gather*}}

The package can be downloaded  from the web:  sites.google.com/site/bedratyuklp/.

\end{document}